\documentclass{amsart}
\usepackage{amscd,amssymb,graphicx}
\usepackage[all,cmtip,line]{xy}

\theoremstyle{plain}
\newtheorem{thm}{Theorem}[section]
\newtheorem*{thmA}{Theorem 1}
\newtheorem*{thmB}{Theorem 2}
\newtheorem*{prpC}{Proposition 3}
\newtheorem*{thmD}{Theorem 4}
\newtheorem*{thmE}{Theorem 5}
\newtheorem*{thmF}{Theorem 6}

\newtheorem{lem}[thm]{Lemma}
\newtheorem{cor}[thm]{Corollary}
\theoremstyle{remark}
\newtheorem{rem}[thm]{Remark}
\newtheorem{quest}[thm]{Question}
\newcommand{\m}{\phantom{-}}

\newcommand{\N}{\mathbb{N}}
\newcommand{\Z}{\mathbb{Z}}
\newcommand{\R}{\mathbb{R}}
\newcommand{\C}{\mathbb{C}\mkern1mu}
\renewcommand{\H}{\mathbb{H}\mkern1mu}
\newcommand{\Ca}{\mathbb{O}\mkern 1mu}

\newcommand{\Sph}{\mathbb{S}}

\DeclareMathOperator{\diag}{diag}

\DeclareMathOperator{\id}{id}

\newcommand{\E}{{\mathchoice
{\mathrm 1\mskip-4.2mu\mathrm l}{\mathrm 1\mskip-4.2mu\mathrm l}
{\mathrm 1\mskip-3.9mu\mathrm l}{\mathrm 1\mskip-4.0mu\mathrm l}}}

\DeclareMathOperator{\Diff}{Diff}

\renewcommand{\mod}{\mathrm{~mod~}}
\newcommand{\bmat}{\left(\begin{smallmatrix}}
\newcommand{\emat}{\end{smallmatrix}\right)}
\newcommand{\bsmat}{\bigl[\begin{smallmatrix}}
\newcommand{\esmat}{\end{smallmatrix}\bigr]}
\newcommand{\Bsmat}{\Bigl[\begin{smallmatrix}}
\newcommand{\Esmat}{\end{smallmatrix}\Bigr]}
\newcommand{\bbsmat}{\biggl[\begin{smallmatrix}}
\newcommand{\eesmat}{\end{smallmatrix}\biggr]}
\newcommand{\BBsmat}{\Biggl[\begin{smallmatrix}}
\newcommand{\EEsmat}{\end{smallmatrix}\Biggr]}
\newcommand{\mml}{\langle\mkern-3mu\langle}
\newcommand{\mmr}{\rangle\mkern-3mu\rangle}

\newcommand{\abs}[1]{\vert #1\vert}

\newcommand{\SO}{\mathrm{SO}}
\newcommand{\OO}{\mathrm{O}}

\newcommand{\Syp}{\mathrm{Sp}}
\newcommand{\Gtwo}{\mathrm{G}_2}

\newcommand{\GL}{\mathrm{GL}}

\DeclareMathOperator{\Real}{Re}
\DeclareMathOperator{\Imag}{Im}

\newcommand{\GM}{\Sigma^7_{\mathrm{GM}}}
\newcommand{\ip}{\langle \,\cdot\,,\,\cdot\,\rangle}
\begin{document}
\title{An infinite family of Gromoll-Meyer spheres}
\author{Carlos Dur\'an}
\address{IMECC-UNICAMP, Pra\c{c}a Sergio Buarque de Holanda, 651, 
Cidade Universit\'aria - Bar\~ao Geraldo, 
Caixa Postal: 6065 
13083-859 Campinas, SP, Brasil }
\email{cduran@ime.unicamp.br}
\author{Thomas P\smash{\"u}ttmann}
\address{Mathematisches Institut\\Universit\"at Bonn\\
        D-53115 Bonn\\Germany}
\email{puttmann@math.uni-bonn.de}
\author{A.~Rigas}
\address{IMECC-UNICAMP, Pra\c{c}a Sergio Buarque de Holanda, 651, 
Cidade Universit\'aria - Bar\~ao Geraldo, 
Caixa Postal: 6065 
13083-859 Campinas, SP, Brasil }
\email{rigas@ime.unicamp.br}
\thanks{C.~Duran and A.~Rigas were supported by CNPq. C.~Duran was also
supported by FAPESP grant 03/016789 and FAEPEX grant 15406. T.~P\"uttmann
was supported by a DFG Heisenberg fellowship and by the DFG priority program
SPP~1154 ``Globale Differentialgeometrie''.}

\begin{abstract}
We construct a new infinite family of models of exotic $7$-spheres.
These models are direct generalizations of the Gromoll-Meyer sphere.
From their symmetries, geodesics and submanifolds
half of them are closer to the standard $7$-sphere than any other
known model for an exotic $7$-sphere.
\end{abstract}

%\subjclass[2000]{Primary 57T20}

\maketitle

%
% ################################
%
\section{Introduction}
This paper provides a new geometric way to construct all exotic $7$-spheres.
Exotic spheres are differentiable manifolds that are homeomorphic but
not diffeomorphic to standard spheres. The first examples were found by
Milnor \cite{milnor} in 1956 among the $\Sph^3$-bundles over $\Sph^4$.
It turned out that $7$ is the smallest dimension where exotic spheres
can occur except possibly in the special dimension~$4$.
In any dimension $n > 4$ the exotic spheres and the standard sphere
form a finite abelian group: the group $\Theta_n$ of (orientation preserving
diffeomorphism classes of) homotopy spheres \cite{kervaire}.
The inverse element in $\Theta_n$ can be obtained by a
change of orientation. In dimension\,7 we have $\Theta_7\approx \Z_{28}$.
Hence, ignoring orientation there are $14$ exotic $7$-spheres.
From these $14$ exotic $7$-spheres four
(corresponding to $2,5,9,12,16,19,23,26\in\Z_{28}$)
are not diffeomorphic to an $\Sph^3$-bundle over $\Sph^4$ \cite{eells}.

In 1974 Gromoll and Meyer \cite{meyer} constructed an exotic $7$-sphere,
$\GM$, as quotient of the compact group $\Syp(2)$ by a two-sided $\Sph^3$-action.
This construction provided $\GM$ automatically with a metric of
nonnegative sectional curvature ($K \ge 0$). The Gromoll-Meyer sphere $\GM$
was the only exotic sphere known to admit such a metric until 1999
when Grove and Ziller \cite{grove} constructed metrics with $K\ge 0$
on all Milnor spheres, i.e., on all exotic $7$-spheres that are
$\Sph^3$-bundles over $\Sph^4$.
In 2002 Totaro \cite{totaro} and independently Kapovitch and Ziller
\cite{kapovitch} showed that $\GM$ is the only exotic sphere that
can be modeled by a biquotient of a compact group and thus
underlined the singular status of the Gromoll-Meyer sphere among
all models for exotic spheres.

We nevertheless provide an elementary and direct generalization of
the Gromoll-Meyer construction. The essential components in this
construction are natural self-maps of $\Sph^7$, namely, the $n$-powers
of unit octonions, $n\in \Z$. In terms of quaternions these maps are defined by
\begin{gather*}
  \rho_n: \Sph^7 \to \Sph^7, \quad
  \bmat \cos t + p \sin t\\ w \sin t\emat \mapsto \bmat \cos nt + p \sin nt \\ w \sin nt\emat
\end{gather*}
where $p\in \Imag \H$ and $w\in \H$ with $\abs{p}^2 + \abs{w}^2 = 1$.
Let $\mml u,v\mmr := \bar u^{\mathrm{t}} v$ denote the standard
Hermitian product on $\H^2$. The submanifolds
\begin{gather*}
  E_n^{10} := \bigl\{ (u,v) \in \Sph^7\times \Sph^7\;\big\vert\;
    \mml \rho_n(u), v\mmr = 0 \bigr\}
\end{gather*}
come equipped with a free action of the unit quaternions:
\begin{gather*}
  \Sph^3 \times E_n^{10} \to E_n^{10},\quad
  q \star (u,v) = (q u \bar q, q v).
\end{gather*}
Here, $qu\bar q$ means that the two quaternionic components of $u$ are
simultaneously conjugated by $q \in \Sph^3$.
The quotient of $E_n^{10}$ by the free $\star$-action is a smooth manifold
\begin{gather*}
  \Sigma^7_n := E_n^{10}/\Sph^3.
\end{gather*}
For $n = 1$ we have $E_1^{10} = \Syp(2)$ (the group of quaternionic
$2\times 2$ matrices $A$ with $\bar A^{\mathrm{t}} A = \E$) and the
$\star$-action is the original Gromoll-Meyer action. Hence, $\Sigma^7_1 =\GM$.
It is also easy to see that $\Sigma^7_0$ is diffeomorphic to $\Sph^7$.

\begin{thmA}
The differentiable manifold $\Sigma^7_n$ is a homotopy sphere and represents
the $(n\mod 28)$-th element in $\Theta_7\approx\Z_{28}$.
\end{thmA}

Let $\Z_2 \times \Z_2$ denote the diagonal matrices of $\OO(2) \subset \Syp(2)$.
All $E_n^{10}$ admit a smooth action of $\Z_2 \times \Z_2 \times \Sph^3$
that commutes with the free $\star$-action:
\begin{gather*}
  \Z_2 \times \Z_2 \times E_n^{10} \to E_n^{10}, \quad B \bullet (u,v) = (B u, B v),\\
  \Sph^3 \times E_n^{10} \to E_n^{10}, \quad q \bullet (u,v) = (u,v\bar q).
\end{gather*}
The induced effective action on $\Sigma^7_n$ is an action of
$\Z_2\times\Z_2\times\SO(3)$ where $\SO(3) = \Sph^3/\{\pm 1\}$.
On $\Sigma^7_0$ this action can be identified with the linear action
\begin{gather*}
  (B,\pm q)\cdot (x,u) = (Bx,Bqu\bar q)
\end{gather*}
on $\Sph^7 \subset \R^2 \times (\Imag \H)^2$. On $\Sigma^7_1=\GM$ the
action coincides with the subaction of the $\OO(2)\times\SO(3)$-action
given in \cite{meyer}.

\smallskip

The surprising fact is the following even/odd grading of the $\Sigma^7_n$:

\begin{thmB}
All $\Sigma^7_n$ with even $n$ are equivariantly homeomorphic to $\Sph^7$
with the linear $\Z_2\times\Z_2\times\SO(3)$-action given above.
All $\Sigma^7_n$ with odd $n$ are equivariantly homeomorphic to the
Gromoll-Meyer sphere $\GM$ with the above $\Z_2\times\Z_2\times\SO(3)$-action.
If $n$ is even all fixed point sets in $\Sigma^7_n$ are spheres while if
$n$ is odd there are also $3$-dimensional fixed point sets with fundamental
groups $\Z_2$ and $\Z_3$.
\end{thmB}

The even/odd grading of the $\Sigma^7_n$ also transfers to
some of the invariant submanifolds. The most important one is $\Sigma^5_n$
whose preimage under the map $E_n^{10}\to \Sigma^7_n$ consists of points
$(u,v)$ where both quaternionic components of $u$ are purely imaginary.

\begin{prpC}
\label{tc}
$\Sigma^5_n$ is $\Z_2\times\Z_2\times\SO(3)$-equivariantly diffeomorphic to
$\Sph^5\subset (\Imag\H)^2$ with the linear action $(B,\pm q)\cdot u = Bqu\bar q$
if $n$ is even and to the Brieskorn sphere $W^5_3$ if $n$ is odd.
The subsphere $\Sigma^5_n$ is minimal for every
$\{\pm\E\} \times \SO(3)$-invariant metric on $\Sigma^7_n$.
\end{prpC}
Recall here that the Brieskorn sphere $W^5_d$ with $d \in \N$ is the intersection
of the unit sphere in $\C^4 =\C \oplus \C^3$ with the complex hypersurface
\begin{gather*}
  z_0^3 + z_1^2 + z_2^2 + z_3^2 \;=\; 0
\end{gather*}
and that there is a natural $\OO(2) \times \SO(3)$-action on~$W^5_d$:
\begin{align}
\begin{split}
\label{actionbries}
  \OO(2) \times \SO(3) \times W^5_d &\to W^5_d, \\
  \bigl( \bmat \cos \theta & -\sin \theta\\ \sin\theta & \m \cos\theta \emat, A\bigr)
    \cdot (z_0, z) \;&=\; (e^{2i\theta}z_0, e^{d i\theta} A z),\\
  \bigl( \bmat 1 & \m 0\\ 0 & -1\emat , A\bigr)
    \cdot (z_0, z) \;&=\; ( \bar z_0, A \bar z).
\end{split}
\end{align}
The classification theorems of J\"anich and Hsiang-Hsiang imply that
for $G = \OO(2)\times\SO(3)$ and even for the smaller group
$G = \{\pm \E\} \times\SO(3)$ the Brieskorn sphere $W^5_d$ is not
$G$-equivariantly homeomorphic to $\Sph^5$ with any linear action, see \cite{mayer}.
However, $W^5_d$ is $\SO(3)$-equivariantly diffeomorphic to $\Sph^5$.
In the case $d = 3$ an explicit formula for such a diffeomorphism is given in \cite{dp}.

\smallskip

The invariant subsphere $\Sigma^5_n$ is dual to the invariant circle
$\Sigma^1_n$ whose preimage under the map $E_n^{10}\to \Sigma^7_n$
consists of points $(u,v)$ for which both components of $u$ are real.
These two dual submanifolds play a central role for the geodesic geometry
of $\Sigma^7_n$.
We construct a one parameter family of $\Z_2\times\Z_2\times\SO(3)$-invariant
metrics $\ip_{\nu}$ on each $\Sigma^7_n$ with the following property:

\begin{thmD}
All points $p\in \Sigma^1_n$ have the wiedersehen property, i.e.,
every unit speed geodesic $\gamma$ in $\Sigma^7_n$ with
$\gamma(0) = p$ is length minimizing on $[0,\pi[$ and obeys
$\gamma(\pi) = -p$ and $\gamma(2\pi) = p$.
Moreover, $\Sigma^1_n$ and $\Sigma^5_n$ have constant
distance $\tfrac{\pi}{2}$ and the map
$\Sigma^1_n \ast \Sigma^5_n \to \Sigma^7_n$
that maps $(x,y,t)$ to $\gamma(t)$, where $\gamma: [0,\tfrac{\pi}{2}] \to \Sigma^7$
is the unique unit speed geodesic segment from $x$ to $y$, is a homeomorphism.
\end{thmD}

This invariant geodesic join structure actually is the key to prove
Theorem\,1 and Theorem\,2. In the particular case of exotic $7$-spheres 
our method is an improvement over the general construction that equips
all exotic spheres with pointed wiedersehen metrics \cite{besse}.

\medskip

The even/odd grading of the $\Sigma^7_n$ is in contrast to what happens
for the Milnor spheres $M^7_{k,l}$ and the Brieskorn spheres $W^7_{6n-1,3}$.

The Milnor sphere $M^7_{k,l}$ with $k+l = 1$ is defined by gluing two
copies of $\H \times \Sph^3$ along $(\H\smallsetminus \{0\})\times\Sph^3$
by the map
\begin{gather}
\label{gluemilnor}
(u,v) \mapsto
  \bigl(\tfrac{u}{\abs{u}^2}, (\tfrac{u}{\abs{u}})^k \,v \, (\tfrac{u}{\abs{u}})^l \bigr).
\end{gather}
For convenience, we set $M^7_{k,l} = M^7_d$ where $d = k-l$ is odd.
The Milnor sphere $M^7_d$ represents the $\tfrac{d^2-1}{8}$-th
element in $\Theta_7$, see \cite{eells}.
There is a natural $\SO(3) = \Sph^3/\{\pm 1\}$-action on $M^7_d$
which is in both charts defined~by
\begin{gather*}
  \pm q \bullet (u,v) = (qu\bar q, q v \bar q).
\end{gather*}
Davis \cite{davis} has shown that $M^7_d$ is $\SO(3)$-equivariantly
diffeomorphic to $M^7_{d'}$ if and only if $d' = \pm d$ and that all $M^7_d$
are $\SO(3)$-equivariantly homeomorphic to $\Sph^7\subset\H^2$
with the linear $\SO(3)$-action given by $(\pm q,u) \mapsto qu\bar q$.
We show that the latter situation changes when one extends the
$\SO(3)$-action by the commuting involution
\begin{gather*}
   (u,v) \mapsto (u,-v)
\end{gather*}
in both charts. This involution fixes all points in the base of the bundle
$M^7_d\to \Sph^4$ and induces the antipodal map on all the $\Sph^3$-fibers.
For consistency, the group generated by $\SO(3)$ and the involution
is denoted by $\{\pm \E\}\times \SO(3)$.

\begin{thmE}
The fixed point set of the involution $(-\E,\pm i)$ on $M^7_d$
is a $3$-dimensional lens space with fundamental group $\Z_{\abs{d}}$.
Hence, $M^7_d$ is $\{\pm \E\}\times \SO(3)$-equivariantly homeomorphic
to $M^7_{d'}$ if and only if $d = \pm d'$. Moreover, for $\abs{d} > 3$ none of
the $M^7_d$ are $\{\pm \E\}\times\SO(3)$-equivariantly homeomorphic to any
of the $\Sigma^7_n$.
\end{thmE}

This theorem is a consequence of Theorem\,\ref{milnorbries}
which is the analogue of Proposition\,3 for the Milnor spheres.

Grove and Ziller \cite{grove} constructed $\SO(3)$-actions on $M^7_d$
that are entirely different from the $\SO(3)$-actions on $M^7_d$ and
$\Sigma^7_n$ above. The $\SO(3)$-actions on $M^7_d$ and $\Sigma^7_n$
fix a circle pointwise while the Grove-Ziller actions are almost free.

\smallskip

The Brieskorn sphere $W^7_{6n-1,3}$ is defined by the intersection of
the the unit sphere $\Sph^9\subset \C^5 = \C\oplus\C\oplus \C^3$
with the complex hypersurface
\begin{gather*}
  w^{6n-1}+z_0^3 + z_1^2 + z_2^2 + z_3^2 = 0.
\end{gather*}
It represents the $(n\mod 28)$-th homotopy sphere in $\Theta_7$
(see \cite{brieskorn}) and admits the natural $\SO(3)$-action
$\bigl(A,(w,z_0,z)\bigr) \mapsto (w,z_0,Az)$.

\begin{thmF}
None of the $W^7_{6n-1,3}$ are $\SO(3)$-equivariantly
diffeomorphic to any of the~$\Sigma^7_n$ or to any of the $M^7_{k,l}$.
\end{thmF}

In particular, $W^7_{6n-1,3}$ is not $\SO(3)$-equivariantly homeomorphic
to the join of a circle and $W^5_3$. Thus, the equivariant topology of the
$\Sigma^7_n$ with odd $n$ is much more determined by the equivariant topology
of $W^5_3$ than the equivariant topology of $W^7_{6n-1,3}$
although the latter contain $W^5_3$ in a much more obvious way
(just by setting $w = 0$).

\smallskip

Many of the constructions in this paper generalize the constructions
given in \cite{dp} for the original Gromoll-Meyer sphere $\GM$.

\medskip

The authors would like to thank Uwe Abresch for several useful discussions
and Wolfgang Ziller for many valuable suggestions.

\bigskip

\section{A construction of the $\Sph^3$-principal bundles over $\Sph^7$}
\label{ggm}
Recall from the introduction the definition of $E_n^{10}\subset \Sph^7 \times \Sph^7$:
\begin{gather*}
  E_n^{10} = \bigl\{ (u,v) \in \Sph^7\times \Sph^7\;\big\vert\;
    \mml \rho_n(u), v\mmr = 0 \bigr\}.
\end{gather*}
For $n = 1$ the space $E_1^{10}$ can be equivalently seen as the group $\Syp(2)$
of $2\times 2$ quaternionic matrices $A$ such that $\bar A^{\mathrm{t}} A = \E$.
The standard projection $\Syp(2) \to \Sph^7$, $A = (u,v) \mapsto u$ turns
$\Syp(2)$ into an $\Sph^3$-principal bundle over $\Sph^7$.

\begin{lem}
$E_n^{10}$ is the pull-back of $\Syp(2)$ by the map $\rho_n: \Sph^7\to \Sph^n$.
\end{lem}
\begin{proof}
By the usual explicit construction, the total space of the pull-back bundle
$\rho_n^{\ast}\bigl(\Syp(2)\bigr)$ is the submanifold of $\Sph^7\times \Syp(2)$
consisting of all pairs $(u,A)$ such that $\rho_n(u)$ is the first colum of~$A$.
It is evident, however, that in this construction we log the first column of $A$ twice.
Eliminating this redundancy leads to the definition of $E_n^{10}$ above.
This in particular shows that $E_n^{10}$ is a submanifold of $\Sph^7\times \Sph^7$.
\end{proof}

\begin{cor}
$E_n^{10}$ is an $\Sph^3$-principal bundle over $\Sph^7$
classified by $n\mod 12$.
\end{cor}

\begin{proof}
The $\Sph^3$-principal bundles over $\Sph^7$ are classified by
$\pi_6(\Sph^3)\approx\Z_{12}$ and the characteristic map of the
bundle $\Syp(2)\to \Sph^7$ generates $\pi_6(\Sph^3)$ (see
\cite{hu} or \cite{duranrigas} for a more explicit reference).
The map $\rho_n$ has degree $n$.
\end{proof}

The principal bundle map $E_n^{10}\to \Sph^7$ is given by the projection
to the first column. The corresponding free $\Sph^3$-action on $E_n^{10}$
is given by
\begin{gather*}
  \Sph^3 \times E_n^{10} \to E_n^{10},\quad q\bullet (u,v) = (u,v\bar q).
\end{gather*}
The map $\tilde\rho_n$ in the pull-back diagram
\begin{gather*}
\begin{CD}
E_n^{10} @>{\tilde\rho_n}>> \Syp(2) \\
@VVV @VVV \\
\Sph^7 @>{\rho_n}>> \Sph^7
\end{CD}
\end{gather*}
takes the explicit form
\begin{gather*}
  \tilde\rho_n : E_n^{10} \to \Syp(2), \quad (u,v) \mapsto (\rho_n(u),v).
\end{gather*}
Recall from the introduction that there is a free $\Sph^3$-action
$q\star (u,v) = (qu\bar q,qv)$ on $E_n^{10}$ that commutes with the
$\bullet$-action and whose orbit space is the smooth manifold $\Sigma^7_n$.
The pull-back diagram above extends to the following commutative diagram:
\begin{gather*}
\xymatrix{
  & E_n^{10} \ar[rr]^{\tilde \rho_n} \ar[dr] \ar[dd] & & \Syp(2) \ar[dr] \ar[dd] & \\
  & & \Sigma^7_n \ar[rr]^{\rho_n'\qquad} & & \GM\\
  & \Sph^7 \ar[rr]^{\rho_n} & & \Sph^7 & }
\end{gather*}
The degree of the induced map $\rho_n' : \Sigma^7_n \to \GM$ is $n$.
The proof that $\Sigma^7_n$ represents the $(n \mod 28)$-th element
of $\Theta_7$ requires several geometric constructions and is postponed
until section\,\ref{geodjoin}.

\smallskip

Each principal bundle $E_n^{10}$ admits a natural action of
$\Z_2 \times \Z_2 \times \Sph^3\times\Sph^3\times\Sph^3$,
where $\Z_2 \times \Z_2$ denotes the diagonal matrices in $\OO(2) \subset \Syp(2)$:
\begin{gather}
\label{fullactiona}
  \Z_2 \times \Z_2 \times E_n^{10} \to E_n^{10}, \quad B \cdot (u,v) = (B u, B v),\\
\label{fullactionb}
  \Sph^3\times\Sph^3\times\Sph^3 \times E_n^{10} \to E_n^{10}, \quad
  (q_1,q_2,q_3) \cdot \bmat u_1 & v_1\\ u_2 & v_2 \emat
    = \bmat q_1 u_1 \bar q_1 & q_1 v_1 \bar q_3\\
      q_2 u_2 \bar q_1 & q_2 v_2 \bar q_3\emat.
\end{gather}

\begin{lem}
This action on $E_n^{10}$ is of cohomogeneity~2.
\end{lem}
\begin{proof}
The third $\Sph^3$-factor yields the principal action related to the
bundle $E_n^{10} \to \Sph^7$, $(u,v) \mapsto u$, i.e.,
this $\Sph^3$-factor acts simply transitively on the fiber over any $u\in \Sph^7$.
The action of the first two $\Sph^3$-factors on $\Sph^7$ has kernel
$\{\pm (1,1)\}$ and induces a standard linear $\SO(4)$-action on $\Sph^7$.
By applying all three $\Sph^3$-factors one can transform an
arbitrary point in $E_n^{10}$ to a point of the form
\begin{gather*}
  \bmat \cos t + i \cos s \sin t & -\sin s \sin nt\\ \sin s\sin t & \cos nt - i \cos s \sin nt\emat.
  \qedhere
\end{gather*}
\end{proof}

The diagonal in the first two $\Sph^3$-factors gives the Gromoll-Meyer
action $\star$ corresponding to the principal bundle $E_n^{10}\to\Sigma^7_n$.
The third $\Sph^3$-factor and the $\Z_2\times \Z_2$-factor yield the effective
$\Z_2\times\Z_2\times\SO(3)$-action $\bullet$ on~$\Sigma^7_n$ from the introduction.
It is an interesting question for which $n$ this $\bullet$-action can be extended.
The maximum dimension of any compact differentiable transformation
group of an exotic $7$-sphere is four \cite{straume}.
On the original Gromoll-Meyer sphere $\GM = \Sigma^7_1$
there is a natural $\OO(2)\times \SO(3)$-action. It is induced by the action
\begin{gather*}
  \OO(2) \times \SO(3) \times \Syp(2) \to \Syp(2), \quad
   (A,q) \bullet (u,v) \mapsto (Au,Av\bar q)
\end{gather*}
on $\Syp(2) = E_1^{10}$ and extends the $\bullet$-action naturally.
A corresponding $\OO(2)\times \SO(3)$-action exists of course on $\Sigma^7_{-1}$.
On $\Sigma^7_0$ an $\OO(2)\times \SO(3)$-action is induced by
the action
\begin{gather*}
\OO(2) \times \SO(3) \times E_0^{10} \to E_0^{10}, \quad
   (A,q) \bullet (u,v) \mapsto (Au, v\bar q).
\end{gather*}
On the other $\Sigma^7_n$ with $n \neq -1,0,1$, however, it seems likely that the
$\Z_2\times\Z_2\times \SO(3)$-action cannot be extended to any larger group,
see Remark\,\ref{iso}.

\begin{quest}
Which $E_n^{10}$ admit Riemannian metrics with $K\ge 0$ that are
invariant under the cohomogeneity~2 action above? If some $E_n^{10}$
admits such a metric then by the O'Neill formulas the induced metric on
$\Sigma^7_n$ also has $K\ge 0$. This would be particularly interesting for
those $\Sigma^7_n$ that are not diffeomorphic to $\Sph^3$-bundles over
$\Sph^4$ since on such exotic spheres no metrics with $K \ge 0$ are known so far.
\end{quest}

\begin{rem}
While there are twelve $\Sph^3$-principal bundles over $\Sph^7$ there
are $28$ homotopy $7$-spheres. This means in particular that some
$\Sigma^7_n$ are quotients of trivial bundles $E_n^{10}$.
This phenomenon is well-known from surgery theory (see \cite{wall})
in an inexplicit~way.
\end{rem}

\begin{rem}
Grove and Ziller \cite{grove} constructed cohomogeneity one metrics
with $K \ge 0$ on all $\Sph^3\times\Sph^3$-principal bundles over $\Sph^4$.
It is known that the $(n \mod 12)$-th $\Sph^3$-principal bundle over $\Sph^4$
is diffeomorphic to an $\Sph^3\times\Sph^3$-principal bundles over $\Sph^4$
if and only if $n \mod 12 \in \{0,1,3,4,6,7,9,10\}$. It is easy to see that the set of
all integers $n$ with $n \mod 12 \in \{0,1,3,4,6,7,9,10\}$ maps
surjectively on $\Z_{28}$. Thus, every element in $\Theta_7$ can
be represented by some $\Sigma^7_n$ such that $E_n^{10}$
admits a cohomogeneity one metric with $K \ge 0$. However,
this does not mean that $\Sigma^7_n$ admits a metric with $K\ge 0$
since the Gromoll-Meyer action $E_n^{10}$ is not isometric with respect to
the Grove-Ziller metric.
\end{rem}

\bigskip

\section{Invariant submanifolds and parity}
In this section we will see that the even/odd grading of the generalized
Gromoll-Meyer spheres $\Sigma^7_n$ is based on an elementary property
of the maps~$\rho_n$.

Consider the subsets
\begin{align*}
 E^9_n &:= \{(u,v) \in E_n^{10} \;\vert\; u \in \Imag \H \times \H\},\\
 E^8_n &:= \{(u,v) \in E_n^{10} \;\vert\; u \in \Imag \H \times \Imag \H\}
\end{align*}
of $E_n^{10} \subset \Sph^7\times \Sph^7$. These are the preimages of the subspheres
\begin{align*}
  \Sph^6 &= \bigl\{ \bmat p\\ w \emat \;\vert\;
    p\in \Imag \H,~w \in \H, ~\abs{p}^2 + \abs{w}^2 = 1 \bigr\},\\
  \Sph^5 &= \bigl\{ \bmat p_1\\ p_2 \emat \;\vert\;
    p_1,p_2\in \Imag \H,~\abs{p_1}^2 + \abs{p_2}^2 = 1 \bigr\}
\end{align*}
of $\Sph^7 \subset \H\times \H$ under the principal
bundle projection $E_n^{10} \to \Sph^7$.

\begin{lem}
$E^9_n$ and $E^8_n$ are submanifolds of $E_n^{10}$ diffeomorphic to
$\Sph^6\times\Sph^3$ and $\Sph^5\times \Sph^3$, respectively.
\end{lem}

\begin{proof}
$E_n^9\to \Sph^6$ is a proper subbundle of $E_n^{10} \to \Sph^7$ and hence trivial.
\end{proof}

\begin{lem}
$E_n^9$ and $E_n^8$ are invariant under the free $\star$-action of $\Sph^3$ and
under the $\bullet$-action of $\Z_2\times\Z_2\times \Sph^3$.
Hence, the $\star$-quotients $\Sigma^6_n$ and $\Sigma^5_n$ are submanifolds
of $\Sigma^7_n$ with a natural $\bullet$-action of $\Z_2\times\Z_2\times \SO(3)$.
\end{lem}
\begin{proof}
Straightforward.
\end{proof}

\begin{lem}
As submanifolds of $\Sph^7\times\Sph^7$ we have
\begin{align*}
\ldots = E^9_{-3} = E^9_{-1} &= E^9_1 = E^9_3 = \ldots,\\
\ldots = E^9_{-4} = E^9_{-2} = & E^9_0 = E^9_2 = E^9_4 = \ldots
\end{align*}
and the same identities also hold for $E^8_n \subset E^9_n$ and
for the quotients $\Sigma^6_n$ and~$\Sigma^5_n$.
\end{lem}

\begin{proof}
This is an immediate consequence of the two basic identities
\begin{gather*}
  \rho_{2m+1}\bigl(\bmat p\\ w\emat\bigr) =  (-1)^m \bmat p\\ w\emat
  \quad\text{and}\quad
  \rho_{2m}\bigl(\bmat p\\ w\emat\bigr) =  (-1)^m \bmat 1\\ 0\emat
\end{gather*}
for $\bmat p\\ w\emat \in \Sph^6 \subset \Imag\H \times \H$.
\end{proof}

\begin{cor}
\label{fivesphere}
If $n$ is odd, $\Sigma^5_n$ is equivariantly diffeomorphic to the Brieskorn
sphere $W^5_3$ with its natural $\Z_2\times\Z_2\times \SO(3)$-action.
If $n$ is even, $\Sigma^5_n$ is equivariantly diffeomorphic to the Euclidean
sphere $\Sph^5 \subset \R^3\times \R^3$ where $\SO(3)$-acts diagonally on
both $\R^3$-factors and each $\Z_2$-factor acts on one of the $\R^3$-factors.
\end{cor}

\begin{proof}
From Lemma\,7.4 of \cite{dp} we know that $\Sigma^5_1$
is $\Z_2\times\Z_2\times \SO(3)$-equivariantly diffeomorphic to
 $W^5_3$. For $\Sigma^5_0$ we observe that
\begin{gather*}
E^8_0 = \bigl\{ \bmat p_1 & 0 \\ p_2 & q \emat 
  \;\vert\; p_1, p_2 \in \Imag\H, q \in \Sph^3 \bigr\}
\end{gather*}
The natural embedding $\Sph^5\to E^8_0$, $\bmat p_1\\ p_2\emat
\mapsto \bmat p_1 & 0 \\ p_2 & 1 \emat$
identifies the $\star$-quotient of $E^8_0$ with~$\Sph^5$.
\end{proof}

\begin{lem}
The subsphere $\Sigma^5_n$ is minimal in $\Sigma^6_n$ and
$\Sigma^7_n$ for all $\{\pm \E\} \times \SO(3)$-invariant
Riemannian metrics on $\Sigma^7_n$.
\end{lem}

\begin{proof}
Analogous to the proof of Corollary\,3.4 in \cite{dp} this follows from
the fact that $\Sigma^5_n$ is the union of orbits whose isotropy
groups contain elements of the form~$(-\E,\pm q)$.
\end{proof}

\bigskip

\section{The geodesic join structure of $\Sigma^7_n$}
\label{geodjoin}
We will now study the geometry of a one parameter family of Riemannian
metrics on $E_n^{10}$ and $\Sigma^7_n$ and use the results to prove
Theorem\,1, Theorem\,2 and Theorem\,4. The one parameter family of
metrics is defined in such a way that the constructions of \cite{duran}
and \cite{dp} for $\GM$ can be extended to all $\Sigma^7_n$.

We equip the total space of the principal bundle $E_n^{10} \to \Sph^7$
with the Riemannian metric $\ip_{\nu}$ with $\nu > 0$
defined by the following properties:
\begin{itemize}
\item
  The $\Sph^3$-fibers have constant curvature $\tfrac{1}{\nu}$.
\item
  The horizontal distribution is given by the pull-back of
  the horizontal distribution of $\Syp(2)$ via the map $\rho_n$,
  i.e., we pull-back the principal bundle connection of $\Syp(2)$.
\item
  The metric $\ip_{\nu}$ induces on $\Sph^7$ the metric with
  constant curvature~$1$ by Riemannian submersion.
\end{itemize}
Such metrics are called connection metrics or Kaluza-Klein metrics.

The $\Z_2\times \Z_2\times \Sph^3\times \Sph^3\times\Sph^3$-action
given in (\ref{fullactiona}) and (\ref{fullactionb}) is isometric with respect
to the metric $\ip_{\nu}$. In particular, the Gromoll-Meyer action $\star$ is isometric
and $\Sigma^7_n$ inherits a Riemannian metric by Riemannian submersion,
which will again be denoted by $\ip_{\nu}$. The $\bullet$-action of
$\Z_2\times \Z_2\times \Sph^3$ on $E_n^{10}$ is also isometric.
Since the $\bullet$-action commutes with the $\star$-action, it induces an
effective isometric $\Z_2\times\Z_2\times \SO(3)$-action on $(\Sigma^7_n,\ip_{\nu})$.

\begin{lem}
The common fixed point set of $\SO(3)$ in $\Sigma^7_n$ is the circle
\begin{gather*}
  \Sigma^1_n := \bigl\{\pi_{\Sigma^7_n}
    \bigl( \bmat \cos t & -\sin nt\\ \sin t & \m \cos nt \emat \bigr)
    \;\vert\; t \in \R \bigr\}.
\end{gather*}
Hence, for any $\SO(3)$-invariant Riemannian metric on $\Sigma^7_n$,
this circle $\Sigma^1_n$ is a simple closed geodesic.
\end{lem}

\begin{proof}
$\pi_{\Sigma^7_n}(u,v)$ is a fixed point of $\SO(3)$ if and only if
for every $q \in \Sph^3$ there is a $q' \in \Sph^3$ such that
$(q' u \bar q', q' v \bar q) = (u,v)$. It is easy to see from the second
column of this equation that all elements of $\Sph^3$ occur for $q'$.
Therefore, $u$ must have two real components.
\end{proof}
Note that the $\Z_2\times\Z_2$-action on $\Sigma^1_n$ is equivalent
to the standard $\Z_2\times\Z_2$-action on $\Sph^1$. In particular,
for each point $p\in \Sigma^1_n$ there is a natural antipode $-p$.

\begin{thm}
\label{wiedersehen}
Every unit speed geodesic $\gamma$ in $(\Sigma^7_n,\ip_{\nu})$
with $\gamma(0) = p \in \Sigma^1_n$ is length minimizing
on $[0,\pi[$ and we have $\gamma(\pi) = -p$ and $\gamma(2\pi) = p$.
\end{thm}

\begin{proof}
The proof is similar to the proofs of Theorem\,I in \cite{duran} and
Theorem\,2.1 in~\cite{dp}.
We lift $\gamma$ horizontally to a geodesic $\tilde \gamma$ in $E_n^{10}$ with
\begin{gather*}
  \tilde\gamma(0) = \alpha(t) := \bmat \cos t & -\sin nt\\ \sin t & \m \cos nt \emat \in E_n^1.
\end{gather*}
That $\tilde\gamma$ is horizontal with respect to $E_n^{10}\to \Sigma^7_n$
means that the geodesic $\tilde\gamma$ passes perpendicularly through
all $\star$-orbits. It is straightforward to check that
\begin{gather*}
  \Sph^3\star \alpha(t) = \Sph^3\bullet \alpha(t).
\end{gather*}
Thus, $\tilde \gamma$ passes perpendicularly through
$\Sph^3\bullet\tilde\gamma(0)$.
A geodesic that passes perpendicularly through one orbit passes
perpendicularly through all orbits. Hence, $\tilde \gamma$ passes
perpendicularly through all $\Sph^3$-orbits of the $\bullet$-action. 
In other words, $\tilde \gamma$ is horizontal to the principal fibration
$E_n^{10}\to \Sph^7$.
Hence, $\tilde\gamma$ projects to a geodesic $\beta$ in $\Sph^7$.
By definition of $\ip_{\nu}$ the sphere $\Sph^7$ inherits the metric with constant
curvature $1$ from $E_n^{10}$ by Riemannian submersion.
Since all unit speed geodesics of $\Sph^7$ that start at
$\beta(0) = \pi_{\Sph^7}(\alpha(t))$ pass through $\beta(\pi)=-\beta(0)$
at time $\pi$ we have $\beta(\pi) = \pi_{\Sph^7}(\alpha(t+\pi))$.
Thus, $\tilde\gamma(\pi)$ is contained in
$\Sph^3\bullet \alpha(t+\pi) = \Sph^3\star \alpha(t+\pi)$ and
$\tilde\gamma(2\pi)$ is contained in
$\Sph^3\star \alpha(t+2\pi) = \Sph^3\bullet \tilde\gamma(0)$.
This shows $\gamma(\pi) = -\gamma(0)$ and $\gamma(2\pi) = \gamma(0)$.
Now let $\gamma$ be a unit speed geodesic in $\Sigma^7_n$ with $\gamma(0) = p$
and $\gamma_1(l) = -p$. By the construction above $\beta$ is a unit speed
geodesic in $\Sph^7$ with $\beta(l)= -\beta(0)$. Hence, $l$ cannot be less than $\pi$.
\end{proof}

Recall that the join $X\ast Y$ of two spaces $X$ and $Y$ is the quotient
of $X\times Y\times [0,1] / \sim$ where $(x,y,0) \sim (x,y',0)$ and
$(x,y,1) \sim (x',y,1)$ for all $x\in X$ and all $y\in Y$.
For our purposes it is convenient to substitute $[0,1]$ by $[0,\tfrac{\pi}{2}]$.
\begin{cor}
\label{joingeod}
$\Sigma^1_n$ and $\Sigma^5_n$ have constant distance $\tfrac{\pi}{2}$.
Moreover, the map $\Sigma^1_n \ast \Sigma^5_n \to \Sigma^7_n$
that maps $(x,y,t)$ to $\gamma(t)$, where $\gamma: [0,\tfrac{\pi}{2}] \to \Sigma^7$
is the unique unit speed geodesic segment from $x$ to $y$, is an equivariant
homeomorphism.
\end{cor}

\begin{proof}
This follows from the construction in the proof of Theorem\,\ref{wiedersehen}
if one recalls that the submanifolds $E_n^1$ and $E_n^9$ of $E_n^{10}$
project to the submanifolds
\begin{align*}
  \Sph^1 &= \bigl\{ \bmat \cos t\\ \sin t\emat \;\big\vert\; t \in \R \bigr\},\\
  \Sph^5 &= \bigl\{ \bmat p_1\\ p_2 \emat \;\big\vert\;
    p_1,p_2\in \Imag \H,~\abs{p_1}^2 + \abs{p_2}^2 = 1 \bigr\}
\end{align*}
of $\Sph^7\subset \H^2$ under the principal fibration $E_n^{10} \to \Sph^7$
and to the submanifolds $\Sigma^1_n$ and $\Sigma^5_n$ of $\Sigma^7_n$
under the principal fibration $E_n^{10} \to \Sigma^7_n$.
\end{proof}

Theorem\,\ref{wiedersehen} and Corollary\,\ref{joingeod} together
yield Theorem\,4 from the introduction.

\begin{cor}
\label{eqjoin}
$\Sigma^7_n$ is $\Z_2\times\Z_2\times \SO(3)$-equivariantly homeomorphic
to $\Sph^1 \ast \Sph^5$ if $n$ is even and to $\Sph^1 \ast W^5_3$ if $n$ is odd.
Here, the $\Z_2\times\Z_2$ acts on $\Sph^1$ in the standard way.
\end{cor}
\begin{proof}
This is evident from Corollary\,\ref{fivesphere} and Corollary\,\ref{joingeod}.
\end{proof}

In particular, all $\Sigma^7_n$ with even $n$ are mutually equivariantly
homeomorphic and that all $\Sigma^7_n$ with odd $n$ are mutually
equivariantly homeomorphic. This proves Theorem\,2 from the introduction.

\medskip

\begin{proof}[Proof of Theorem\,1]
Consider the unit speed geodesic
\begin{gather*}
 \beta(t) =
   \bmat \cos t + p\sin t \\ w \sin t \emat
\end{gather*}
in $\Sph^7 \subset \H^2$ that emanates from the north pole with initial velocity
$\bmat p\\ w \emat\in \Sph^6\subset \Imag\H \times \H$.
A lift $\tilde \gamma_n$ of this curve to $E_n^{10}$ with
$\tilde \gamma_n(0) = \bmat 1 & 0\\ 0 & 1\emat$ is given by
\begin{gather}
\label{lift}
  \tilde \gamma_n (t) =
    \bmat \cos t + p\sin t &\; -e^{ntp} \bar w \sin(nt)\\
    w \sin t &\;  
    {\tfrac{w}{\abs{w}} e^{ntp}(\cos(nt) - p \sin(nt)) \tfrac{\bar w}{\abs{w}}} \emat.
\end{gather}
Here $e^{p} = \cos \abs{p} + \tfrac{p}{\abs{p}} \sin \abs{p}$
denotes the exponential map of $\Sph^3\subset\H$ at $1$.
Note that for $w = 0$ equation (\ref{lift}) simply becomes
$\tilde\gamma_n(t) = \bmat e^{ntp} & 0\\ 0 & 1\emat$.
Using the identity
\begin{gather*}
  \tilde\rho_n(\tilde \gamma_n(t)) = \tilde\gamma_1(nt)
\end{gather*}
for the map $\tilde\rho_n : E_n^{10} \to \Syp(2)$ defined in section\,\ref{ggm}
it is straightforward to verify that $\tilde \gamma_n$ is the unique {\em horizontal}
lift of $\beta$ to $E_n^{10}$ with $\tilde \gamma_n(0) = \E$.
Since the fibers of $E_n^{10}\to \Sph^7$ and $E_n^{10}\to \Sigma^7_n$
through $\tilde\gamma_n(0)  = \bmat 1 & 0\\ 0 & 1\emat$ are the same (as sets),
the geodesic $\tilde\gamma_n$ is horizontal with respect to both these fibrations.
This  shows that $\gamma_n = \pi_{\Sigma^7_n}\circ \tilde \gamma_1$ is a geodesic
in $\Sigma^7_n$.
Now, considering all possible unit initial vectors
$\bmat p\\ w \emat\in \Sph^6\subset \Imag\H \times \H$ and times
$t \in [0,\tfrac{\pi}{2}]$ the geodesics $\gamma_n$ provide an embedding of a
disk $D^7(\tfrac{\pi}{2})$ into $\Sigma^7_n$ by Theorem\,\ref{wiedersehen}.
In the same way, the geodesics
$\pi_{\Sigma^7}\circ(- \tilde\gamma_n)\circ (-\id)$
provide another embedding of the same disk. By Theorem\,\ref{wiedersehen},
$\Sigma^7_n$ is the twisted sphere obtained by gluing these two embedded
disks along their common boundary. One easily checks that
\begin{gather*}
  \tilde\gamma_n(p,w,\tfrac{\pi}{2})
  = q\star (-\tilde\gamma_n(- p',-w',\tfrac{\pi}{2}))
\end{gather*}
for some $q\in\Sph^3$ if and only if  $(p',w') = \sigma^n(p,w)$
where $\sigma$ is the exotic diffeomorphism of
$\Sph^6 \subset \Imag\H\times \H$ first described in \cite{duran}.
This diffeomorphism $\sigma$ generates $\pi_0(\Diff_+(\Sph^6))$.
It is given by the formula
\begin{gather*}
  \sigma(p,w) :=  \overline{ \mathrm{b}(p,w)} (p,w) \mathrm{b}(p,w)
\end{gather*}
where $\mathrm{b}(p,w) = \tfrac{w}{\abs{w}} e^{\pi p} \tfrac{\bar w}{\abs{w}}$
is an analytic formula for a generator of $\pi_6(\Sph^3)$. Hence, we have
obtained $\Sigma^7_n$ by gluing two $7$-disks with the $n$-th power of
a generator of $\pi_0(\Diff_+(\Sph^6)) \approx \Theta_7 \approx \Z_{28}$.
\end{proof}

\begin{rem}
\label{iso}
Let $G$ be a compact group acting smoothly on $\Sigma^7_n$ with
$\Z_2\times\Z_2\times \SO(3) \subset G$. Precisely as in \cite{dp}, Lemma\,3.7,
one can show that $G$ leaves $\Sigma^1_n$ and $\Sigma^5_n$ invariant.
Let $n\not\in \{-1,0,1\}$. Comparing for different $p\in\Sigma^1_n$ the closing
behaviour of geodesics that start at $p$ perpendicularly to $\Sigma^1_n$,
one can see that $\Z_2\times\Z_2$ is the maximal compact group that acts
isometrically on $(\Sigma^7_n,\ip_{\nu})$ and effectively on the circle $\Sigma^1_n$.
This difference from the cases $n=-1,0,1$ suggests that
$\Z_2\times\Z_2\times \SO(3)$ is the full isometry group of $(\Sigma^7_n,\ip_{\nu})$.
\end{rem}

\begin{rem}
If we pull back the metric $\ip_{\nu}$ on $\Syp(2)$ by the map $\tilde\rho_n$
then we obtain a degenerate metric $\ip_{\nu}'$ on $E_n^{10}$ that has
the same geodesics through the circle $\Sigma^1_n$ as the metric $\ip_{\nu}$.
For $n \not \in \{-1,0,1\}$ the metric $\ip_{\nu}'$ is degenerate precisely over
$\abs{n}-1$ subspheres in $\Sph^7$ whose first quaternionic
components have constant real part. With such a metric $\Sigma^7_n$
looks like $n$ copies of $\GM$ stacked one on top of the other,
i.e., like a degenerate connected sum of $n$ copies of $\GM$.
\end{rem}

\begin{rem}
The manifolds $(E^9_n, \ip_{\nu})$ with even $n$ are not just mutually equal
as submanifolds of $\Sph^7\times \Sph^7$ but also mutually equal as
Riemannian manifolds. Hence, also the manifolds $(\Sigma^6_n, \ip_{\nu})$
with even $n$ are all mutually equal as Riemannian manifolds.
The analogous statements hold for odd $n$.
\end{rem}

\bigskip

\section{Comparison to the exotic Milnor and Brieskorn 7-spheres}
\label{comp}
In this section we compare the equivariant topology of the spheres $\Sigma^7_n$
with the equivariant topology of the Milnor spheres $M^7_d$ and the
Brieskorn spheres $W^7_{6n-1,3}$ and prove Theorem\,5 and Theorem\,6
of the introduction.

\smallskip

Recall from the introduction that the Milnor spheres $M^7_d$ admit
natural $\{\pm\E\} \times \SO(3)$-actions. Davis \cite{davis} has shown that
these actions can be extendend to $\GL(2,\R)\times \SO(3)$-actions.
In the first chart the $\GL(2,\R)$-action is given by
\begin{gather*}
  \bmat a & c\\ b & d\emat \bullet (u,v)
    = \Bigl(\tfrac{au+c}{bu+d}, \det\bmat a & c\\ b & d\emat 
       \bigl(\tfrac{bu+d}{\abs{bu+d}}\bigr)^k v  \,
       \bigl(\tfrac{bu+d}{\abs{bu+d}}\bigr)^{l} \Bigr)
\end{gather*}
and in the second one by
\begin{gather*}
  \bmat a & c\\ b & d\emat \bullet (u,v)
    = \Bigl(\tfrac{b+du}{a+cu}, \det\bmat a & c\\ b & d\emat 
       \bigl(\tfrac{a+c\bar u}{\abs{a+c\bar u}}\bigr)^k v \,
       \bigl(\tfrac{a+c\bar u}{\abs{a+c\bar u}}\bigr)^{l} \Bigr).
\end{gather*}
Note that our definition of the action differs from the definition
given by Davis by the factor $\det\bmat a & c\\ b & d\emat$.
The reason is that with our definition the identification between
$M^7_3$ and the Gromoll-Meyer sphere $\GM$ given in \cite{meyer}
becomes an $\OO(2)\times\SO(3)$-equivariant diffeomorphism
while without the determiant factor the identification is only
$\SO(2)\times\SO(3)$-equivariant. Moreover, note that the
map $M^7_d \to M^7_{-d}$ given by $(u,v) \mapsto (\bar u,\bar v)$
in both charts is an $\GL(2,\R)\times\SO(3)$-equivariant diffeomorphism.

\begin{thm}
\label{milnorbries}
In every Milnor sphere $M^7_d$ there is a unique invariant submanifold
$M^5_d$ which is $\OO(2)\times\SO(3)$-equivariantly diffeomorphic to the
Brieskorn sphere~$W^5_{\abs d}$ with the $\OO(2)\times\SO(3)$-action
given in (\ref{actionbries}). This submanifold $M^5_d$ is minimal for any
$\{\pm\E\}\times\SO(3)$-invariant Riemannian metric on $M^7_d$.
\end{thm}

\begin{proof}
It suffices to consider the case $d > 0$. Let $M^5_d$ be the submanifold
of $M^7_d$ given by the equations $\Real v = 0$ and $\Real uv = 0$ 
n both charts (it is essential here that $k+l = 1$).
Hirsch and Milnor \cite{hirsch} proved that $M^5_d$ is homeomorphic
and hence (because exotic spheres do not exist in dimension $5$)
diffeomorphic to $\Sph^5$. It is straightforward to check that
$M^5_d$ is invariant under the $\SO(2)\times \SO(3)$-action.
Consider the curve $\alpha$ in $M^5_d$ which is given by
$\alpha(s) = (i \tan s,j)$ in the first chart.
The isotropy groups along $\alpha$ are
\begin{gather*}
   K_{-} =  \{(\E, \pm e^{j\tau})\} \cup \{(-\E, \pm i e^{j\tau})\}
   \cup \bigl\{\bigl( \bsmat 1 & \m 0\\ 0 & -1\esmat, \pm e^{j\tau}\bigr)\bigr\}
   \cup \bigl\{\bigl( \bsmat -1 & \m 0\\ \m 0 & \m 1\esmat, \pm i e^{j\tau}\bigr)\bigr\}
\end{gather*}
at $s = 0$,
\begin{gather*}
  H = \bigl\{(\E,\pm 1), (-\E,\pm i), \bigl( \bsmat 1 & \m 0\\ 0 & -1\esmat, \pm j\bigr),
  \bigl( \bsmat -1 & \m 0\\ \m 0 & \m 1\esmat, \pm k\bigr) \bigr\}
\end{gather*}
for $0 < s < \tfrac{\pi}{4}$, and
\begin{gather*}
  K_{+} = \bigl\{ \bigl( \bsmat \cos \theta & -\sin \theta\\ \sin \theta & \m \cos \theta\esmat,
    \pm e^{-\frac{d}{2}i\theta} \bigr)\bigr\}
    \cup  \bigl\{ \bigl( \bsmat \cos \theta & -\sin \theta\\ \sin \theta &
    \m \cos \theta\esmat \cdot \bsmat 1 & \m 0\\ 0 & -1\esmat,
    \pm e^{-\frac{d}{2}i\theta}j \bigr)\bigr\}
\end{gather*}
at $s = \tfrac{\pi}{4}$. Now consider the Brieskorn sphere $W^5_d$ with
the $\OO(2)\times \SO(3)$-action given in (\ref{actionbries}) and the curve
\begin{gather*}
  \beta(s) = \Bigl( s,0,\tfrac{1}{\sqrt{2}} \sqrt{1-s^2-s^d},-\tfrac{i}{\sqrt{2}}\sqrt{1-s^2+s^d} \Bigr)
\end{gather*}
on the interval $[s_-,0]$ where $s_- < 0$ is the root of $1-s^2+s^d$.
Straightforward computations show that the isotropy
groups along $\beta$ are the same as the isotropy groups along $\alpha$.
This proves that $M^5_d$ and $W^5_d$ are equivariantly diffeomorphic.
The uniqueness and minimality of $M^5_d$ follows from the following fact:
The fixed point set of any element of the form $(-\E,\pm q)$ is contained
in $M^5_d$ and even more $M^5_d$ can be seen
to be the union of orbits whose isotropy groups contains such elements.
\end{proof}

\begin{proof}[Proof of Theorem\,5]
The involution $(-\E,\pm i)$ is contained in $M^5_d \approx W^5_{\abs{d}}$.
The fixed point set of $(-\E,\pm i) = (-\E,\diag(1,-1,-1)$ in $W^5_{\abs{d}}$
is the $W^3_{\abs{d}}$ given by the equation $z_1 = 0$ and hence
diffeomorphic to a lens space with fundamental group $\Z_{\abs{d}}$.
\end{proof}

The Milnor sphere $M^7_d$ have direct analogues $M^{15}_d$ in dimension $15$.
They are obtained by gluing two copies of $\Ca\times \Sph^7$ along
$(\Ca\smallsetminus\{0\})\times \Sph^7$ by the map (\ref{gluemilnor}).
Precisely as above each $M^{15}_d$ admits a smooth action of
$\OO(2)\times \Gtwo$ (see \cite{davis}).
\begin{thm}
In every $M^{15}_d$ there is a unique invariant submanifold $M^{13}_d$
which is $\OO(2)\times\Gtwo$-equivariantly diffeomorphic to the Brieskorn
sphere~$W^{13}_{\abs d}$ with the action of
$\OO(2)\times\Gtwo \subset \OO(2)\times \SO(7)$
given analogously to (\ref{actionbries}).
This submanifold $M^{13}_d$ is minimal for any
$\{\pm\E\}\times\Gtwo$-invariant Riemannian metric on $M^{15}_d$.
\end{thm}
\begin{proof}
Analogous to the proof of Theorem\,\ref{milnorbries}.
\end{proof}

\medskip

Finally, we turn to the Brieskorn spheres $W^7_{6n-1,3}$ and prove
Theorem\,6 from the introduction.
\begin{proof}[Proof of Theorem\,6]
The involution
$\bmat 1 & \m 0 & \m 0\\ 0 & -1 & \m 0\\ 0 & \m 0 & -1\emat \in \SO(3)$
on $W^7_{6n-1,3}$ is given by
\begin{gather*}
  (w,z_0,z_1,z_2,z_3) \mapsto (w,z_0,z_1,-z_2,-z_3).
\end{gather*}
Its fixed point set is thus identical to $W^3_{6n-1,3,2}$, which is the
intersection of the unit sphere $\Sph^5$ in $\C^3$ with the complex hypersurface
\begin{gather*}
  w^{6n-1}+z_0^3 + z_1^2 = 0.
\end{gather*}
Milnor \cite{milnortwo} has shown that $W^3_{5,3,2}$ is
diffeomorphic to Poincare dodecahedral space and that
the universal covering space of $W^3_{6n-1,3,2}$ is
non-compact if $n > 1$.
\end{proof}

\bigskip

%
% ************ bibliography ************
%

\end{document}